\documentclass[12pt]{amsart}
\usepackage{graphicx}
\usepackage{amssymb}
\usepackage{amsfonts}
\usepackage{amsmath}
\usepackage{array}
\usepackage{pstricks}
\usepackage{bm}
\usepackage{epsfig}
\usepackage[utf8]{inputenc}

\usepackage{pspicture}
\usepackage[matrix,ps,xdvi]{xypic} % ne pas enlevera

\newdimen\Squaresize \Squaresize=14pt
\newdimen\Thickness \Thickness=0.5pt

\def\Square#1{\hbox{\vrule width \Thickness
   \vbox to \Squaresize{\hrule height \Thickness\vss
      \hbox to \Squaresize{\hss#1\hss}
   \vss\hrule height\Thickness}
\unskip\vrule width \Thickness}
\kern-\Thickness}

\def\Vsquare#1{\vbox{\Square{$#1$}}\kern-\Thickness}

\def\boxit#1#2{\setbox1=\hbox{\kern#1{#2}\kern#1}%
\dimen1=\ht1 \advance\dimen1 by #1 \dimen2=\dp1 \advance\dimen2 by #1
\setbox1=\hbox{\vrule height\dimen1 depth\dimen2\box1\vrule}%
\setbox1=\vbox{\hrule\box1\hrule}%
\advance\dimen1 by .4pt \ht1=\dimen1
\advance\dimen2 by .4pt \dp1=\dimen2 \box1\relax}

\newdimen\vcadre\vcadre=0.1cm % marges verticales de la boite
\newdimen\hcadre\hcadre=0.1cm % marges horizontales de la boite
\def\arx#1[#2]{\ifcase#1 \relax \or%
  \ar @{-}[#2]  \or%
  \ar @2{-}[#2] \or%
  \ar @{--}[#2] \or%
  \ar @2{.}[#2] \or%
  \ar @{~}[#2]  \fi}

\newtheorem{example}{Example}[section]

\newtheorem{theorem}[example]{Theorem}

\newtheorem{corollary}[example]{Corollary}

\newtheorem{proposition}[example]{Proposition}

\newtheorem{lemma}[example]{Lemma}

\def\Proof{\noindent \it Proof -- \rm}
\def\qed{\hspace{3.5mm} \hfill \vbox{\hrule height 3pt depth 2 pt width 2mm}
\bigskip}

\def\FQSym{{\bf FQSym}}

\def\Sym{{\bf Sym}}

\def\ssh{\Cup}

\def\<{\langle}
\def\>{\rangle}

\def\K{\operatorname{\mathbb K}}

\def\F{{\bf F}}

\def\G{{\bf G}}

\def\SG{{\mathfrak S}}

\def\CC{{\bf C}}

\def\Des{\operatorname{Des}}

\def\shuff#1#2{\mathbin{
\hbox{\vbox{ \hbox{\vrule \hskip#2 \vrule height#1 width 0pt
}%
\hrule}%
\vbox{ \hbox{\vrule \hskip#2 \vrule height#1 width 0pt
\vrule }%
\hrule}%
}}}

\def\shuf{{\mathchoice{\shuff{7pt}{3.5pt}}%
{\shuff{6pt}{3pt}}%
{\shuff{4pt}{2pt}}%
{\shuff{3pt}{1.5pt}}}}%
\def\shuffle{\,\shuf\,}

\def\SC{{\rm SC}}
\def\RC{{\rm RC}}

\def\maj{{\rm maj}}

\def\inv{{\rm inv}}

\title[]%
{A noncommutative cycle index \\
and new bases of quasi-symmetric functions\\
and noncommutative symmetric functions}
\author[J.-C.~Novelli, J.-Y.~Thibon, F. Toumazet]%
{Jean-Christophe Novelli, Jean-Yves Thibon, and Fr\'ed\'eric Toumazet}

\address[]{[Novelli, Thibon, Toumazet] Laboratoire d'informatique Gaspard-Monge\\
Universit\'e Paris-Est Marne-la-Vall\'ee \\
5, Boulevard Descartes \\ Champs-sur-Marne \\
77454 Marne-la-Vall\'ee cedex 2 \\
France}
\email[Jean-Christophe Novelli]{novelli@u-pem.fr}
\email[Jean-Yves Thibon]{jyt@u-pem.fr} 
\email[Fr\'ed\'eric Toumazet]{frederic.toumazet@u-pem.fr} 
\date{\today}


\date{\today}

\begin{document}

\begin{abstract}
We define a new basis of the algebra of quasi-symmetric functions by lifting
the cycle-index polynomials of symmetric groups to  noncommutative polynomials
with coefficients in the algebra of free quasi-symmetric functions, and then
projecting the coefficients to $QSym$.  By duality, we obtain a basis of
noncommutative symmetric functions, for which a product formula and a
recurrence in the form of a combinatorial complex are obtained. This basis
allows to identify noncommutative symmetric functions with the quotient of
$\FQSym$ induced by the pattern-replacement relation $321\equiv 231$ and
$312\equiv 132$.
\end{abstract}

\maketitle
%{\footnotesize
%\tableofcontents
%}
%%%%%%%%%%%%%%%%%%%%%%%%%%%%%%%%%%%%%%%%%%%%%%%%%%%%%%%%%%%%%%%%%%%%%%%%%%%%%%%
%%%%%%%%%%%%%%%%%%%%%%%%%%%%%%%%%%%%%%%%%%%%%%%%%%%%%%%%%%%%%%%%%%%%%%%%%%%%%%%
%%%%%%%%%%%%%%%%%%%%%%%%%%%%%%%%%%%%%%%%%%%%%%%%%%%%%%%%%%%%%%%%%%%%%%%%%%%%%%%
\section{Introduction}

For a partition $\mu=1^{m_1}2^{m_2}\cdots n^{m_n}$ of $n$, the coefficient
of $p_\mu=p_1^{m_1}p_2^{m_2}\cdots p_n^{m_n}$ in the polynomial
$h_n(p_1,\ldots,p_n)$ defined by $h_0=1$ and the recursion 
\begin{equation}
nh_n = h_{n-1}p_1+h_{n-2}p_2+\cdots+p_n
\end{equation}
is the number of permutations of cycle type $\mu$, divided by $n!$. Thus,
$h_n$ is the cycle index of the symmetric group $\SG_n$, and if the $p_k$
are interpreted as the power sums symmetric functions, then $h_n$ is the complete
symmetric function (sum of all monomials of degree $n$), {\it cf.} \cite{Mcd}.

There is a noncommutative analogue of this polynomial, which defines the
so-called noncommutative power sums of the first kind $\Psi_n$  in terms of
the noncommutative complete functions $S_n$ \cite{NCSF1}
\begin{equation}
nS_{n} = S_{n-1}\Psi_1+S_{n-2}\Psi_2+\cdots+\Psi_n.
\end{equation}
It is not difficult to see that the coefficient $c_I$ in
\begin{equation}
  n! S_n
= \sum_I c_I\Psi^I
= \sum_I \frac{n!}{i_1(i_1+i_2)\cdots(i_1+i_2+\cdots i_r)}\Psi^I
\end{equation}
counts the permutations with {\it ordered} cycle type $I$, that is, the
composition consisting of the lengths of the cycles, ordered by increasing
values of their maxima.

For example,
\begin{align}
4!h_4
 &= p_{1111}+6p_{211}+3p_{22}+8p_{31}+6p_4,\\
4!S_4
 &=\Psi^{1111}+3\Psi^{112}+2\Psi^{121}+6\Psi^{13}+\Psi^{211}+3\Psi^{22}+2\Psi^{31}+6\Psi_4,
\end{align}
the 3 transpositions of ordered cycle type $(1,1,2)$ being $(1)(2)(43)$,
$(1)(3)(42)$ and $(2)(3)(41)$, while $(1)(32)(4)$ and $(2)(31)(4)$ are of type
$(1,2,1)$, and $(21)(3)(4)$ is of type $(2,1,1)$.

Setting $y_n=(n-1)!p_n$ and $Y_n=(n-1)!\Psi_n$, which amounts to forgetting
the order inside the cycles, $n!h_n$ and $n!S_n$ become respectively the usual
exponential Bell polynomial $B_n(y_1,\ldots,y_n)$, whose coefficients count
set partitions according to the sizes of the blocks, and its natural
noncommutative analogue \cite{SR} refining this counting according to a
canonical ordering of the blocks.

This situation has been exploited is \cite{NTT} in order to define further
generalizations of the Bell polynomials, whose coefficients live in the
algebra of free quasi-symmetric functions and encode the set partitions
counted by the original coefficients in a way preserving their algebraic
properties. The commutative images of these coefficients are quasi-symmetric
functions, which turn out to form the dual immaculate basis of \cite{BBSSZ} up
to reversal of the compositions.

The aim of this note is to apply the same strategy to the noncommutative cycle
index with its original variables, that is, we keep $Y_n=\Psi_n$, and define
successively a $q$-analogue $c_I(q)$ of $c_I$, a free quasi-symmetric function
$\CC_I$, a quasi-symmetric function $C_I$, and finally, a basis
$U_I=C_{\overline{I}^\sim}$ of $QSym$, and its dual basis $V_I$ for
noncommutative symmetric functions.

This paper is a continuation of \cite{NTT}, to which the reader may refer for
background and notation.

%%%%%%%%%%%%%%%%%%%%%%%%%%%%%%%%%%%%%%%%%%%%%%%%%%%%%%%%%%%%%%%%%%%%%%%%%%%%%%%
%%%%%%%%%%%%%%%%%%%%%%%%%%%%%%%%%%%%%%%%%%%%%%%%%%%%%%%%%%%%%%%%%%%%%%%%%%%%%%%
%%%%%%%%%%%%%%%%%%%%%%%%%%%%%%%%%%%%%%%%%%%%%%%%%%%%%%%%%%%%%%%%%%%%%%%%%%%%%%%
\section{The saillances statistics, and a $q$-analogue}

As the ordered cycle type appears to be of little significance from a group
theoretical point of view, it can be advantageously replaced by a natural
statistic having the same distribution.

Following \cite[Chap. 10]{Loth}, let us say that a word $w=a_1a_2\cdots a_n$
over the integers is {\it initially dominated} if $a_1>a_2,\ldots,a_n$. Order
these words with respect to their maximal letter.  Then, a permutation  has a
unique increasing factorisation
\begin{equation}
	\sigma=u_1u_2\cdots u_r
\end{equation}
into initially dominated words. The composition
\begin{equation}
	\SC(\sigma)=(|u_1|,|u_2|,\ldots, |u_r|)
\end{equation}
will be called the {\it saillance composition} of $\sigma$.

For example, for $\sigma = 351274698$, the factorisation is
$3\cdot 512\cdot  746 \cdot 98$, and $\SC(\sigma)=(1,3,3,2)$.

The natural $q$-analogue of $\Psi_n$ is \cite{NCSF2}
\begin{equation}
	\Theta_n(q)=\frac{S_n((1-q)A)}{1-q}=\sum_{k=0}^{n-1}(-q)^kR_{1^k,n-k} 
\end{equation}
and the expression
\begin{equation}
	[n]_q! S_n(A)   
= \sum_Ic_I(q)\Theta^I(q)
= \sum_I q^{\maj(I)}\frac{[n]_q!}{[i_1]_q[i_1+i_2]_q\cdots[i_1+i_2+\cdots i_r]_q}\Theta^I(q) 
\end{equation}
provides immediately a $q$-analogue of $c_I$:

\begin{proposition}
The polynomial $c_I(q)$ counts the permutations of saillance composition $I$
by number of non-inversions:
\begin{equation}
c_I(q)=\sum_{\SC(\sigma)=I}q^{\binom{n}{2}-\inv(\sigma)}.
\end{equation}
\end{proposition}

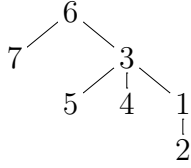
\begin{figure}[ht]
\entrymodifiers={+<4pt>}
\begin{equation*}
\vcenter{\xymatrix@C=4mm@R=2mm{
& {6}\ar@{-}[dl]\ar@{-}[dr] \\
{7} && {3}\ar@{-}[dl]\ar@{-}[d]\ar@{-}[dr] \\
& {5} & {4} & {1}\ar@{-}[d] \\
&&& {2} \\
}}
\end{equation*}
\caption{The poset $P(2,3,2)$.} \label{fig:poset}
\end{figure}

To prove this, we first observe that the inverses of the permutations having a
given saillance composition are the linear extensions of a poset, whose Hasse
diagram is a tree of a very special shape (see Figure \ref{fig:poset}).

Indeed, the conditions $\sigma_1>\sigma_2,\ldots,\sigma_{i_1}$,
$\sigma_1<\sigma_{i_1+1}> \sigma_{i_1+2},\ldots, \sigma_{i_1+i_2}$, etc.,
translate into the conditions that in the inverse permutation, 1 must be to
the right of $2,3,\ldots,i_1$, and that $i_1+1$ must be to the right of
$i_1+2,\ldots, i_1+i_2$ and of $1,2,\ldots i_1$, and so on. This defines a
poset $P(I)$, where $i<_P j$ iff $j$ must be to the right of $i$ in
$\sigma^{-1}$.

For example, for $I=(2,3,2)$, $2<_P 1$; $1<_P3>_P 4,5$; $3<_P6>_P7$, which
gives the comb-tree of Figure \ref{fig:poset}.

Thus, we can apply the Bj\"orner-Wachs $q$-hook-length formula \cite{BW} to
count the linear extensions of this poset by number of inversions.
Since this is clearly a recursive labeling, this gives
\begin{equation}
q^{\inv(\alpha_I)}\frac{[n]_q!}{[i_1]_q[i_1+i_2]_q\cdots[i_1+i_2+\cdots i_r]_q}
\end{equation}
where
\begin{equation}
\alpha_I =  2\cdots i_11\cdot (i_1+2)\cdots (i_1+i_2)\cdot (i_1+1)\cdots 
\end{equation}
is the minimal linear extension, whose number of inversions is
$\sum_k(i_k-1)= n-\ell(I)$.

Thus,
\begin{align}
	\sum_{\SC(\sigma)=I}q^{\inv(\sigma)}&=
	q^{n-\ell(I)}\frac{[n]_q!}{[i_1]_q[i_1+i_2]_q\cdots[i_1+i_2+\cdots i_r]_q}\\
	&=q^{\binom{n}{2}}c_I(q^{-1})
\end{align}
(since $q^{\maj(I)}=q^{i_1}q^{i_1+i_2}\cdots q^{i_1+\cdots+i_{r-1}}$), so that 
\begin{equation}
	c_I(q)= \sum_{\SC(\sigma)=I}q^{\binom{n}{2}-\inv(\sigma)}.
\end{equation}

For $P(2,3,2)$ this polynomial is
\begin{equation}
q^{17} + 3 q^{16} + 6 q^{15} + 9 q^{14} + 11 q^{13} + 12 q^{12} 
+ 11 q^{11} + 9 q^{10} + 6 q^{9} + 3 q^{8} + q^{7}.
\end{equation}

Alternatively, we may replace $\Theta_n(q)$ by $\tilde\Theta_n(q):=(q-1)^{-1}S_n((q-1)A)$.
Then, the coefficient $\tilde c_I(q)$ in
\begin{equation}
	[n]_q! S_n(A)   
= \sum_I\tilde c_I(q)\tilde \Theta^I(q)
= \sum_I \frac{[n]_q!}{[i_1]_q[i_1+i_2]_q\cdots[i_1+i_2+\cdots i_r]_q}\tilde \Theta^I(q) 
\end{equation}
counts permutations of ordered cycle type $I$ by Carlitz's statistic \cite{Car}, \cite[A129178]{Sloane}
\begin{equation}
{\rm invc}(\sigma) = {\rm inv}(f_1(\sigma))
\end{equation}
where $f_1$ is Foata's first fundamental transformation:  express $\sigma$ in standard cycle form 
(i.e., cycles ordered by increasing smallest elements with each cycle written with its smallest element in the first position), 
then remove the parentheses \cite{Loth}.

%%%%%%%%%%%%%%%%%%%%%%%%%%%%%%%%%%%%%%%%%%%%%%%%%%%%%%%%%%%%%%%%%%%%%%%%%%%%%%%
%%%%%%%%%%%%%%%%%%%%%%%%%%%%%%%%%%%%%%%%%%%%%%%%%%%%%%%%%%%%%%%%%%%%%%%%%%%%%%%
%%%%%%%%%%%%%%%%%%%%%%%%%%%%%%%%%%%%%%%%%%%%%%%%%%%%%%%%%%%%%%%%%%%%%%%%%%%%%%%
\section{A dendriform exponential}
%Loday exponential series without denominators
We can now define in $\FQSym$, as in \cite{NTT} for the Bell polynomials,
\begin{equation}
\CC_I := \sum_{\SC(\sigma)=I}\G_\sigma,
\end{equation}
and investigate the commutative image of these elements in $QSym$.

The $\CC_I$ can be expressed in terms of the dendriform operations of
$\FQSym$. We have
\begin{equation}
\CC_n =\G_1\prec \G_1^{n-1}
\end{equation}
and
\begin{equation}
\CC_I = (\cdots (\CC_{i_1}\succ \CC_{i_2})\succ\cdots)\succ \CC_{i_r}
\end{equation}
Their generating series  in $\FQSym\otimes \K\langle Y_1,Y_2,\cdots\rangle$
\begin{equation}
Z := \sum_I \CC_I Y^I 
\end{equation}
satisfies therefore
\begin{equation}
Z = 1 + Z\succ C,\ {\rm with}\ C:=\sum_{n\ge 1}\CC_n Y_n.
\end{equation}

Thus, $\CC_I$ belongs to the Loday-Ronco algebra ${\bf PBT}$: it is the sum of
all  ${\bf P}_T$ where $T$ has a left branch with $r=\ell(I)$ nodes, on which
are grafted as right subtrees all binary trees of sizes
$i_1-1,i_2-1,\ldots,i_r-1$, starting from the bottom.

\begin{figure}
	\footnotesize
%	\begin{mdframed}
\begin{align*}
	Z_1 &= \G_1Y_1\\
	Z_2&= \G_{12}Y^{11}+\G_{21}Y_2\\
	Z_3&=\G_{123}Y^{111}+(\G_{132}+\G_{231})Y^{12} +\G_{213}Y^{21}+ (\G_{312}+\G_{321})Y_3\\
	Z_4&=\G_{1234}Y^{1111}+(\G_{1243} + \G_{1342} +\G_{2341})Y^{112} +(\G_{1324} +\G_{2314})Y^{121} +\G_{2134}Y^{211}\\
	&+(\G_{1423} +\G_{2413} +\G_{3412} +\G_{1432} +\G_{2431} + \G_{3421})Y^{13} +(\G_{2143} +\G_{3241} +\G_{3142})Y^{22}\\
	&+(\G_{4123} +\G_{4132} +\G_{4213} +\G_{4231} +\G_{4312} +  \G_{4321})Y_4 
\end{align*}
%	\end{mdframed}
	\caption{Values of $Z_n$ for $n\le 4$.}\label{fig:Z}
\end{figure}

\begin{proposition}\label{prop:coalg}
The $\CC_I$ span a sub-coalgebra of $\FQSym$, and the coefficients $a_I^{JK}$ in
\begin{equation}\label{eq:coalg}
\Delta \CC_I = \sum_{JK}a_I^{JK}\CC_J\otimes \CC_K
\end{equation}
are nonnegative integers.
\end{proposition}

\Proof
Since $\FQSym$ is self-dual,we can write
\begin{align}
\Delta\CC_I &= \sum_{\alpha,\beta}\<\Delta\CC_I,\F_\alpha\otimes\F_\beta\>\G_\alpha\otimes\G_\beta\\
&=\sum_{\alpha,\beta}\<\CC_I,\F_\alpha\F_\beta\>\G_\alpha\otimes\G_\beta\\
&=\sum_{\SC(\gamma)=I}\sum_{\alpha,\beta}\<\G_\gamma,\F_\alpha\F_\beta\>\G_\alpha\otimes\G_\beta,\\
\end{align}
and it is clear that the distribution of the saillances in the shifted shuffle
$\F_\alpha\F_\beta$ depends only on $\SC(\alpha)$ and $\SC(\beta)$. Thus
\eqref{eq:coalg} holds, and $a_I^{JK}$ is the number of triple of permutations
$(\alpha,\beta,\gamma)$ such that $\SC(\alpha)=J$, $\SC(\beta)=K$,
$\SC(\gamma)=I$, and $\gamma$ occurs in the shifted shuffle of $\alpha$ and
$\beta$.
\qed

A closed formula for $a_I^{JK}$ will be given in the sequel.

Note that $Z$ is an {\it exponential without denominators} according to the
definition of \cite{Lod13} ($Z=e'(C)$ in the notation of this reference), and
that $\sum_n\CC_n=-L(-\G_1)$ is a {\it logarithm without denominators}.

%%%%%%%%%%%%%%%%%%%%%%%%%%%%%%%%%%%%%%%%%%%%%%%%%%%%%%%%%%%%%%%%%%%%%%%%%%%%%%%%
%%%%%%%%%%%%%%%%%%%%%%%%%%%%%%%%%%%%%%%%%%%%%%%%%%%%%%%%%%%%%%%%%%%%%%%%%%%%%%%%
%%%%%%%%%%%%%%%%%%%%%%%%%%%%%%%%%%%%%%%%%%%%%%%%%%%%%%%%%%%%%%%%%%%%%%%%%%%%%%%%
\section{A basis of $QSym$}

We can now take the commutative images of the $\CC_I$, and investigate the
resulting quasi-symmetric functions.

It turns out that $\CC_I(X)$ does not always contain $F_I$, but always
contains $\overline{I}^\sim$. Thus, we define
\begin{equation}
U_I(X) = \sum_{\SC(\sigma)=\overline{I}^\sim}F_{\RC(\sigma)}
\end{equation}
where $\RC(\sigma)$ is the recoil composition of $\sigma$.

The first values are tabulated on Figure \ref{fig:baseU}. This table and all
tables afterwards are quite large and were put at the end of the paper.

\begin{proposition}
The transition matrix $M_n$, whose column indexed by $J$ gives the
coefficients of $U_J$ on the $F_I$ is upper unitriangular when compositions
are ordered lexicographically by length, see Figure \ref{fig:M}.
\end{proposition}

\Proof
Consider a permutation such that $\SC(\sigma)=I$.
Then the values in position $1$, $i_1+1$, $\dots$, $i_1+\dots+i_{r-1}+1$
cannot be recoils since they have no greater value to their left.
Hence, the total number of recoils of $\sigma$ cannot be greater than
$|I|-l(I)$ so that the length of the recoil composition of $\sigma$ cannot be
greater than the length of $\overline{I}^\sim$.
If the length of the recoil composition of $\sigma$ is strictly smaller than
the length of $\overline{I}^\sim$, this composition appears before
$\overline{I}^\sim$ in our order.
Otherwise, all positions except the above mentioned ones are recoils, so that,
given that $\sigma_1$ is greater than $\sigma_k$ for all $k\in[2,i_1]$, the
smallest possible non-recoil of $\sigma$ is at least $i_1$. This means that
all values smaller than $i_1$ are recoils of $\sigma$. Now, if $i_1$ is not a
recoil, the recoil composition of $\sigma$ is strictly smaller than
$\overline{I}^\sim$. Otherwise, they have same first part and one can conclude
by induction on the number of parts of $I$ (the next smallest possible recoil
is then $i_1+i_2$, and so on).
\qed

For example, for $I=(2,3,1,2)$, values at positions $1$, $3$, $6$, and $7$
cannot be recoils so that the total number of recoils cannot be greater than
$8-4=4$, which is the number of recoils of $\overline{I}^\sim=(1,2,1,3,1)$.
If there are exactly 4 recoils, hence at positions $2$, $5$, $4$, $8$,
the smallest non-recoil value is $2$ since $\sigma_1>\sigma_2$, and if $2$ is a
recoil, then the second smallest non-recoil is $5$ since $\sigma_3>\sigma_4$
and $\sigma_3>\sigma_5$. If one computes the largest set of recoils (sorted
by increasing length then by lexicographic order), it is exactly $\{2,5,6,8\}$,
which is $\Des(I)$, so that the recoil composition of any $\sigma$
such that $\SC(\sigma)=I$ is at most (for the specified order on compositions)
its complement $\{1,3,4,7\}=\Des(\overline{I}^\sim)$.

%The minimum polynomial of $M_n$ is $(x-1)^{n\choose 2}$. Indeed, 
%$U_I-F_I$ is a sum of terms of lower major index. Thus, 
%$M_n-I_n$ is the matrix of a nilpotent endomorphism of $QSym$, 
%whose order of nilpotency is at most ${n\choose 2}$. 
%Moreover, $(M_n-I_n)^{{n\choose 2}-1}$ sends $F_{1^n}$ to a nonzero
%multiple of $F_{1,n-1}$.

%%%%%%%%%%%%%%%%%%%%%%%%%%%%%%%%%%%%%%%%%%%%%%%%%%%%%%%%%%%%%%%%%%%%%%%%%%%%%%%%
%%%%%%%%%%%%%%%%%%%%%%%%%%%%%%%%%%%%%%%%%%%%%%%%%%%%%%%%%%%%%%%%%%%%%%%%%%%%%%%%
%%%%%%%%%%%%%%%%%%%%%%%%%%%%%%%%%%%%%%%%%%%%%%%%%%%%%%%%%%%%%%%%%%%%%%%%%%%%%%%%
\section{The dual basis in $\Sym$}

Let $V_I$ denote the dual basis of $U_I$. Reading row $I$ of $M_n$, we have
the expansion of the ribbon $R_I$ on the $V_J$, which is therefore nonnegative.

The first values of the $V_I$ on the ribbon basis are given on Figure
\ref{fig:VR}.

As we shall see, it is rather the expansion on the elementary basis
$\Lambda^I$ which is relevant to the understanding of the $V_I$. The first
values are given on Figure \ref{fig:VL}.

To investigate the multiplicative structure of the $V_I$, we first observe
that $\Sym$ can be identified (as an algebra) with the quotient of $\FQSym$ by
the ideal ${\mathcal J}$ generated by
\begin{equation}
\{\F_\sigma-\F_\tau\,|\, \SC(\sigma)=\SC(\tau)\}.
\end{equation}
Indeed, since the commutative image map $\pi:\ \FQSym\rightarrow QSym$ is an
epimorphism of Hopf algebras, we have 
\begin{equation} 
\Delta\CC_I(X) = \sum_{JK}a_I^{JK}\CC_J(X)\otimes \CC_K(X)
\end{equation}
with the same coefficients $a_I^{JK}$ as in \eqref{eq:coalg}.
Thus,
\begin{equation} 
\Delta U_I(X) = \sum_{JK}a_{\overline{I}^\sim}^{\overline{J}^\sim\overline{K}^\sim} U_J(X)\otimes U_K(X),
\end{equation}
and $V_I$ may be identified with the class
$\overline{\F}_\sigma=\F_\sigma\mod {\mathcal J}$
for $\SC(\sigma)=\overline{I}^\sim$.

Note that this construction is similar to that of Tevlin's fundamental basis
given in \cite{HNTT}.
Moreover, the equivalence classes here can be described in terms of
pattern-replacement relations:

\begin{proposition}\label{prop:congr}
Two permutations $\sigma,\tau\in\SG_n$ have the same saillance composition iff
they are equivalent modulo the pattern replacement relations
\begin{equation}
321\equiv 231\ \text{\ and\ }\ 312\equiv 132.
\end{equation}
\end{proposition}
These relations are the mirror images of those of \cite[Section 3.11]{Kus},
see Section \ref{sec:app} for a proof of the proposition and further details
about this equivalence.

The structure constants
\begin{equation}
V_I V_J=\sum_K c_{IJ}^K V_K
\end{equation}
are therefore given by the following rule. Take two permutations
$\alpha$ with $SC(\alpha)=\overline{I}^\sim$,
such that $\beta$ with $SC(\beta)=\overline{J}^\sim$. Then,
\begin{equation}
c_{IJ}^K=\#\{\F_\gamma\in \F_\alpha\F_\beta| SC(\gamma)=\overline{K}^\sim\}.
\end{equation}

This can be made more explicit. First, we can state a Pieri formula.

For a composition $I$ of $n$, denote by $I[k]$ the composition of $k$ whose
ribbon diagram consists of the first $k$ boxes of that of $I$.

\begin{proposition}\label{prop:pieri}
Let $V'_I=V_{\overline{I}^\sim}$. Then
\begin{equation}
V'_I\Lambda_k = \sum_{j=0}^n \binom{k+j-1}{k-1}V'_{I[n-j],k+j}.
\end{equation}
\end{proposition}

\Proof
For $\sigma\in\SG_n$, the saillance composition of a permutation $\tau$
occuring in the shifted shuffle of $\sigma$ with $\omega_k=k\cdots 21$ is
determined by the position of $m=n+k$. The number of permutations in this
shuffle for which $m$ is at position $n-j$ is $\binom{k+j-1}{k-1}$.
\qed

For example,
\begin{equation}
V'_{32}\Lambda_3 = \binom{2}{2} V'_{323} + \binom{3}{2} V'_{314}
   + \binom{4}{2} V'_{35} + \binom{5}{2} V'_{26} + \binom{6}{2} V'_{17}
   + \binom{7}{2} V'_{8},
\end{equation}
so that, complementing the compositions,
\begin{equation}
V_{1121}V_{111} = 21 \, V_{11111111} + 6 \, V_{1121111} + V_{112211} +
3 \, V_{113111} + 10 \, V_{1211111} + 15 \, V_{2111111}.
\end{equation}

\begin{corollary}
For $I=(i_1,\ldots,i_r)$,
\begin{equation}\label{eq:comp}
V'_I = \sum_{k=i_r}^n(-1)^{k-i_r} \binom{k-1}{i_r-1} V'_{I[n-k]} V'_k.
\end{equation}
\end{corollary}

\Proof If we expand the products on the r.h.s. by the rule of Proposition
\ref{prop:pieri}, the coefficient of $V'_{I[n-p],p}$ is 
\begin{equation}
\sum_{k=i_r}^p(-1)^k \binom{k-1}{i_r-1}\binom{p-1}{k-1}
=(-1)^{i_r}\binom{p-1}{i_r-1}
 \sum_{l=0}^{p-i_r}(-1)^l \binom{p-i_r}{l}
=(-1)^{i_r}\delta_{p,{i_r}}.
\end{equation}
\qed

If we evaluate the alternating sum \eqref{eq:comp} step by step, we can
observe that the partial sums are alternatively positive and negative. This
suggests the existence of a combinatorial complex explaining the formula.

\begin{lemma}Let $u,v$ be two words, and $a$ be a letter. Then,
\begin{equation}
uav = \sum_{u_1u_2=u}(-1)^{|u_2|}u_1\shuffle (a(\overline{u_2}\shuffle v)).
\end{equation}
\end{lemma}

\Proof
It is sufficient to prove the lemma for $u=12\cdots k$, $a=k+1$, and
$v=k+2\cdots n$. In this case, all terms, viewed as elements of $\FQSym$ in
the $\F$-basis, are noncommutative symmetric functions.
The l.h.s. is $R_n$, and for $u_1=1\cdots k- i$, the corresponding term of the
sum is $(-1)^i R_{k-i}R_{1^i,n-k}$. Apart from the last one, which is
$R_{1^k,n-k}$, each term is a sum of two ribbons, and two consecutive terms
have exactly one ribbon in common.
\qed

For example, taking $k=2$, we can write
\begin{equation}
1234 = 12\shuffle 34 - 1\shuffle 3(2\shuffle 4)+3(21\shuffle 4)
\end{equation}
which amounts to the identity
\begin{equation}
R_4=R_2R_2 -R_1R_{12}+R_{112} = (R_4+R_{22})-(R_{22}+R_{112})+R_{112}.
\end{equation}
Applying this to the permutation $2143$, whose saillance composition is $22$,
we obtain
\begin{equation}
2143 = 21\shuffle 43 - 2\shuffle 4(1\shuffle 3)+4(12\shuffle 3)
\end{equation}
which translates into
\begin{equation}
V'_{22}= V'_2V'_2 -2 V'_1V'_3+3 V'_4.
\end{equation}

Refining the argument of the proof of Proposition \ref{prop:pieri} yields the
general product rule:

\begin{theorem}
The coefficients $  \bar c_{IJ}^K$ in the product
\begin{equation}
V'_I V'_J=\sum_K \bar c_{IJ}^K V'_K
\end{equation}
where $\ell(I)=m$, $\ell(J)=p$,  $\ell(K)=q$, and $r=q-p$, are given by
\begin{equation}
\bar c_{IJ}^K =
\begin{cases}
{\displaystyle\prod_{i=1}^p \binom{k_i-1}{j_i-1}}\ \text{if $q=p$,}\\
{\displaystyle\prod_{i=1}^p \binom{k_{i+r}-1}{j_i-1}}
  \ \text{if $r>0$, $(k_1,\ldots;k_{r-1})=(i_1,\ldots,i_{r-1})$
          and $k_r\le i_r$,}\\
0\ \text{otherwise}.
\end{cases}
\end{equation}
\end{theorem}

\Proof
To expand a product $V'_IV'_J$, we have to compute the saillance compositions
$K$ of permutations $\nu$ occuring in the shifted shuffle of permutations
$\sigma$ and $\tau$ of respective saillance compositions $I$ and $J$.

We shall first discuss the number of saillances of $\nu$.
The subword containing the $q\in[p,p+r]$ values of the  saillances of $\nu$
consists of a prefix of the saillances of $\sigma$ followed by all shifted
saillances of $\tau$. So, following the notations of the theorem,
$q\in [p,p+r]$.

Now, if $q=p$, the saillances are exactly those of $\tau$ so this case
corresponds to permutations $\nu$ such that $\nu_1=\tau_1+m$. In that case,
between two  saillances of $\nu$, there must be at least the corresponding
values in $\tau$ \emph{and} some values belonging to $\sigma$, hence the
product of binomial coefficients since these values were shuffled in all
possible ways with values of $\sigma$.

If $r=q-p>0$, the first $r-1$ saillances correspond to the saillances of
$\sigma$, which means in particular that the value $\tau_1+m$ cannot appear
too early, so that, the first $(r-1)$ parts of $K$ have to be equal to the
first $(r-1)$ parts of $I$. The next part of $K$ has to satisfy 
$k_r<i_r$, meaning that either we meet $\tau_1+m$ or the next saillance of
$\sigma$. And now the same explanation as before applies: each saillance value
in $\nu$ coming from $\tau$ must be at least as far apart as these values in
$\tau$, since they were shuffled with values of $\sigma$, whence the product of
binomial coefficients.
\qed

For example, let us compute the coefficient of $V'_{51}$ in
$V'_{2}V'_{31}$.

This amounts to computing the number of permutations in
$21\ssh 3124 = 21 \shuffle 5346$ with saillances at positions $\{1,6\}$.
In that case, the first value has to be $5$ and the last one  has to be $6$, all
remaining values being at any possible place, hence $\binom{4}{2}$: put
$21$ and $34$ in any order in the middle positions.

\medskip
Let us now compute the coefficient of $V'_{2131}$ in
$V'_{21}V'_{121}$.

This amounts to computing the number of permutations in
$213\ssh 1324 = 213 \shuffle 4657$ of saillances at positions $\{1,3,4,7\}$.
Since there are four saillances, the first one has to be the $2$ and the
others are $4$, $6$, and $7$. So we need $2$ in position $1$, then $4$ in
position $3$, then $6$ in position $4$, and $7$ in position $7$.
The value $5$ of $4657$ is in position $5$ or $6$ (between values $6$ and
$7$), value $1$ has to be in position $2$ (no saillance there) and value $3$
in position either $5$ or $6$, hence the binomial coefficient
$\binom{2}{1}$ showing that values $3$ and $5$ were shuffled together at
positions $5$ and $6$.

For example,
\begin{align*}
V'_2V'_{31}&= V'_{132} + 3 \, V'_{141} + V'_{231} + V'_{33} + 3 \, V'_{42} + 6 \, V'_{51}\\
V'_{31}V'_2&= 4 \, V'_{15} + 3 \, V'_{24} + V'_{312} + 2 \, V'_{33} + 5 \, V'_{6}\\
V'_{21}V'_{121}&= V'_{1123} + 2 \, V'_{1132} + 3 \, V'_{1141} + V'_{1222} + 2 \, V'_{1231} + V'_{124} + V'_{1321}\\
& + 2 \, V'_{133} + 3 \, V'_{142} + 4 \, V'_{151} + V'_{21121} + V'_{2122} + 2 \, V'_{2131} + V'_{2221} + V'_{223}\\
& + 2 \, V'_{232} + 3 \, V'_{241} + V'_{322} + 2 \, V'_{331} + V'_{421}
\end{align*}

\medskip
The only known basis which appears to be related to  $V_I$ is the $\Psi$ basis:
\begin{proposition}
The $V$-expansion of $\Psi_n$ is
\begin{equation}
\Psi_n =\sum_{I\vDash n}(-1)^{\ell(I)-1}V_I.
\end{equation}
\end{proposition}

\Proof We know that $\Psi_n=\sum_{k=0}^nR_{1^k,n-k}$, and we have in fact
\begin{equation}
R_{1^k,n-k}=\sum_{\ell(I)=k+1,\ I\vDash n}V_I.
\end{equation}

Indeed, the sum of the  permutations whose recoils are exactly $1,2,\ldots,k$ is
$(k+1)(k\cdots 21\shuffle k+2\cdots n)$, and the saillance compositions of these
permutations are precisely all the compositions of $n$ of length $k+1$.
\qed

%%%%%%%%%%%%%%%%%%%%%%%%%%%%%%%%%%%%%%%%%%%%%%%%%%%%%%%%%%%%%%%%%%%%%%%%%%%%%%%%
%%%%%%%%%%%%%%%%%%%%%%%%%%%%%%%%%%%%%%%%%%%%%%%%%%%%%%%%%%%%%%%%%%%%%%%%%%%%%%%%
%%%%%%%%%%%%%%%%%%%%%%%%%%%%%%%%%%%%%%%%%%%%%%%%%%%%%%%%%%%%%%%%%%%%%%%%%%%%%%%%
\section{Appendix: insertion algorithms for some pattern-replacement
relations}\label{sec:app}

This appendix provides the proof of Proposition \ref{prop:congr}, together
with some supplementary material, answering a question asked by Darij Grinberg
\cite{Gri}: to explain why the cardinalities of the equivalence classes of
two relations considered in \cite{Kus} have the same distribution. One of
these relations turns out to be, up to mirror image of the patterns, the one
inducing equality of the saillance compositions on the inverse permutations.

%%%%%%%%%%%%%%%%%%%%%%%%%%%%%%%%%%%%%%%%%%%%%%%%%%%%%%%%%%%%%%%%%%%%%%%%%%%%%%%%
\subsection{Background}

A pattern-replacement relation is an equivalence relation on permutations
defined by sets of patterns of the same size, two permutations being
equivalent if one is obtained from the other by rearranging the letters
forming a pattern of the set so as to form another one from the same set.
Historically, the first example was defined by the two sets $\{132,312\}$ and
$\{213,231\}$: this is the well-known Knuth equivalence, whose classes are the
fibers of the Robinson-Schensted correspondence.

Actually, the Knuth relations are defined on words over a totally ordered
alphabet, and the quotient of the free monoid by the congruence generated by
these relations is the celebrated plactic monoid, which has been for many
years considered as a unique and singular object.  It was the discovery of
quantum groups and crystal bases which led to the understanding that the
plactic monoid was associated with the root systems of type $A$, and that such
objects existed for other types \cite{LLT,Lit96}.

This was however not the end of the story, as the investigation of the
representation-theoretical meaning of quasi-symmetric functions led to the
discovery of the hypoplactic monoid, which, while related in some way to
quantum groups, does not fit in the previous pattern \cite{NCSF4}.

Finally, the investigation of the product rule of the Loday-Ronco Hopf algebra
of planar binary trees led to the sylvester monoid, and other combinatorial
Hopf algebras provided many new examples, for which no representation
theoretical interpretation is known \cite{HNT}.

All these monoids induce pattern-replacement equivalences when restricted to
permutations. These equivalences are rather special, in that they are induced
by congruences on words which are compatible with standardization and
restrictions to intervals, two properties ensuring that they can be used to
define Hopf algebras \cite{Nz,NT}. Also, in each of these cases, there is an
insertion algorithm analogous to the Robinson-Schensted correspondence.

One may therefore wonder whether there are other pattern-replacement
equivalences, not necessarily coming from such congruences, for which there is
still an insertion algorithm, and some other interesting properties such as a
closed formula for the number of classes, or for the cardinality of a class.

These last two points have been thoroughly investigated in the recent
papers~\cite{Kus,KZ,Va,LPRW,PRW,S12}. In this appendix, we provide insertion
algorithms for two examples from \cite{Kus}. Both can be extended to
bijections by introducing a $Q$-symbol, which turns out to be in both cases
increasing trees of a special shape, thus explaining the curious fact that
both equivalences have the same distribution of the cardinalities of classes.

%%%%%%%%%%%%%%%%%%%%%%%%%%%%%%%%%%%%%%%%%%%%%%%%%%%%%%%%%%%%%%%%%%%%%%%%%%%%%%%%
\subsection{The $\{\{ 312, 321\}, \{ 123, 132 \}\}$-equivalence}

Consider the equivalence $\equiv$ generated by the following relations
\cite[Section 3.1]{Kus}: for $a<b<c$,
\begin{equation}
\label{eq-312321}
\begin{split}
cab &\equiv cba \\
abc &\equiv acb.
\end{split}
\end{equation}
We shall define an algorithm sending a permutation to a poset,
whose Hasse diagram will be a labeled tree of a special shape.
First,  define the \emph{W-chain} of $\sigma=\sigma_1\dots\sigma_n$
as the sequence $S=(s_1,\dots,s_k)$ such that
$s_k$ is the position of whichever is rightmost between $1$ and $n$. Any
other $s_i$ satisfies that $\sigma_1\dots\sigma_{s_{i+1}}$ has $\sigma_{s_i}$
as an extremum (the other one being $\sigma_{s_{i+1}}$).
Note that this definition coincides essentially with the $W_w$ sets of
Definition 3.3 of~\cite{Kus}.

Now, given  a W-chain, represent it as a chain poset with $\sigma_{s_1}$ at the
top and $\sigma_{s_k}$ at the bottom, and place all other values $\sigma_j$
of $\sigma$ as leaves of the topmost element of the chain $\sigma_{s_i}$
(regarded as a linear tree) such that $\sigma_j$ belongs to the interval
between $\sigma_{s_{i}}$ and $\sigma_{s_{i+1}}$. Denote the result by
$P(\sigma)$.

For example, with $\sigma=532498617$, the $W$-chain is $[1,3,5,8]$ and the
corresponding values are $[5,2,9,1]$.
Now, $3$ and $4$ end up as leaves of $5$, and all remaining values as leaves of
$2$.

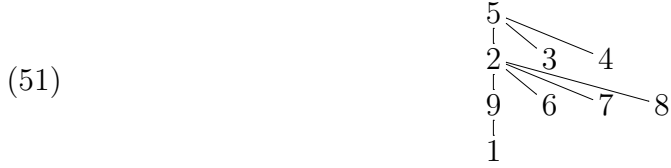
\entrymodifiers={+<4pt>}
\begin{equation}
\vcenter{\xymatrix@C=4mm@R=2mm{
{5}\ar@{-}[d]\ar@{-}[dr]\ar@{-}[drr] \\
{2}\ar@{-}[d]\ar@{-}[dr]\ar@{-}[drr]\ar@{-}[drrr] & {3} & {4} \\
{9}\ar@{-}[d] &  {6} & {7} & {8} \\
{1} \\
      }}
\end{equation}

One can extend this algorithm to a bijection by memorizing in a second
tree of the same shape the position of letter $i$ in $\sigma$.
It follows from the definition of the chain that this yields an
\emph{increasing} tree, that we shall denote by $Q(\sigma)$.

In our example, $Q(\sigma)$ is
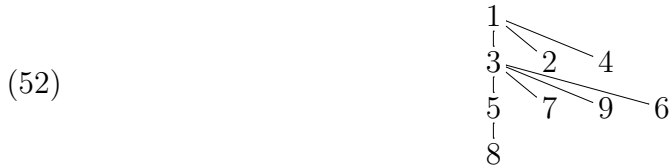
\entrymodifiers={+<4pt>}
\begin{equation}
\vcenter{\xymatrix@C=4mm@R=2mm{
{1}\ar@{-}[d]\ar@{-}[dr]\ar@{-}[drr] \\
{3}\ar@{-}[d]\ar@{-}[dr]\ar@{-}[drr]\ar@{-}[drrr] & {2} & {4} \\
{5}\ar@{-}[d] &  {7} & {9} & {6} \\
{8} \\
      }}
\end{equation}

\begin{theorem}
Consider a permutation $\sigma$.
The linear extensions of $P(\sigma)$, which are in bijection with the
increasing trees of the same shape as $P(\sigma)$, are exactly the
permutations $\equiv$-equivalent to $\sigma$.
\end{theorem}

\Proof
All the necessary material is present in Kuszmaul's
paper~\cite{Kus}: he proves in Lemma~3.4 that if $w\equiv w'$ then they have same $W$-set
hence same poset. And he also proves that if $w$ and $w'$ have same chain,
they have same origin permutation (Lemma 3.6) and that any permutation is
equivalent to its origin permutation (Lemma 3.7).
\qed

Now, given this property, one recovers instantly all the results of \cite{Kus} 
concerning this pattern.

\begin{itemize}
\item There are $2^{n-1}$ classes of permutations of size $n$. Indeed, there are
$2^{n-1}$ possible $W$-sets: given any subset $S$ of $[1,n]$ containing both $1$
and $n$, there are two ways of ordering it as a chain: put either $1$ or $n$
at the bottom and work your way up by taking alternatively the remaining
maximum and minimum in $S$.
\item The size of a class is given by an explicit hook-length formula:
indeed, it is equal to the number of linear extensions of a poset whose Hasse
diagram is a (very special) tree.
\end{itemize}

In our example, the hook-lengths are
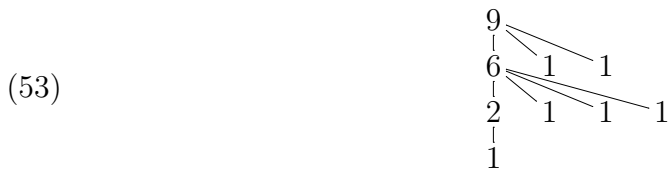
\entrymodifiers={+<4pt>}
\begin{equation}
\vcenter{\xymatrix@C=4mm@R=2mm{
{9}\ar@{-}[d]\ar@{-}[dr]\ar@{-}[drr] \\
{6}\ar@{-}[d]\ar@{-}[dr]\ar@{-}[drr]\ar@{-}[drrr] & {1} & {1} \\
{2}\ar@{-}[d] &  {1} & {1} & {1} \\
{1} \\
      }}
\end{equation}
so that the cardinality of the class is $\frac{9!}{9\cdot 6\cdot 1^6}=3360$.

The counting of equivalence classes can be easily refined as follows. 

\begin{proposition}
Among the $2^{n-1}$ classes,  there are
exactly $\binom{n-1}{k-1}$ classes of permutations beginning with letter $k$. 
\end {proposition}

\Proof
This
amounts to counting (\emph{e.g.}, by induction) $W$-sets beginning with $k$.
More precisely, if one denotes by $A_n^k$ the number of such classes and
${A_n^k}^+$ (resp. ${A_n^k}^-$) the number of such classes whose second
value of the chain is greater (resp. smaller) than the first,
there are exactly $\binom{n-2}{k-1}$ (resp. $\binom{n-2}{k-2}$) such
classes.
\qed

%%%%%%%%%%%%%%%%%%%%%%%%%%%%%%%%%%%%%%%%%%%%%%%%%%%%%%%%%%%%%%%%%%%%%%%%%%%%%%%%
\subsection{The $\{\{ 123,132\}, \{ 213,231 \}\}$-equivalence}

Consider the equivalence $\equiv_2$ generated by the relations
\cite[Section 3.11]{Kus}: for $a<b<c$,
\begin{equation}
\label{eq-231213}
\begin{split}
bac &\equiv_2 bca \\
abc &\equiv_2 acb.
\end{split}
\end{equation}

\def\lrm{\text{lrm}}

Define the \emph{left-to-right minima}, $\lrm$ for short of
$\sigma=\sigma_1\dots\sigma_n$
as the sequence $S=\{s_1,\dots,s_k\}$ such that
$s_i$ is the $i$-th smallest integer such that $\sigma_{s_i}$ has only greater
values to its left in $\sigma$.

\begin{proposition}
The $\equiv_2$ classes are exactly the sets of permutations having the same
left-to-right minima.
\end{proposition}

\Proof
Following~\cite{Kus}, one first checks that in any $\equiv_2$ class, there is
exactly one $V$-permutation (as defined in \cite{Kus}) which is, by the way,
the lexicographically smallest element of the class. This time, the proof is
very easy since one can orient the relations (which is equivalent to saying
that the cardinality of the classes are obtained by pattern avoidance, hence
relating to Theorem 4.7 of~\cite{Kus}) and decide to rewrite any pattern $acb$
into $abc$ and any pattern $bca$ into $bac$. The words having neither $acb$
nor $bca$ patterns are $V$-permutations. This proves that any class has at
least one such element. Conversely, it is obvious that the $\lrm$ of a
permutation does not change with any rewriting, hence there cannot be two
$V$-permutations in the same class.
\qed

Let us now define an algorithm sending a permutation to a poset.
Given the $\lrm$, represent it as a chain poset (again, regarded as a linear
tree) with $\sigma_{s_1}$ at the top and $\sigma_{s_k}$ at its bottom, and
place all other values
$\sigma_j$ of $\sigma$ as leaves of the topmost element of the chain that is
smaller than $\sigma_j$. Denote the result by $P_2(\sigma)$.

For example, with $\sigma=739465281$, its lrm is $[1,2,7,9]$ and the
corresponding values are $[7,3,2,1]$.
Now, $8$ and $9$ end as leaves of $7$, and all remaining values as leaves of
$3$.

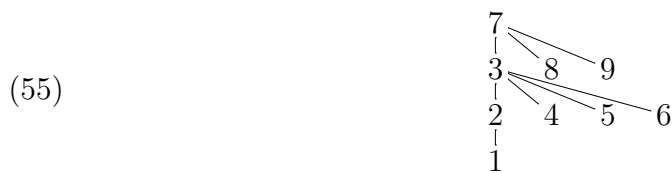
\entrymodifiers={+<4pt>}
\begin{equation}
\vcenter{\xymatrix@C=4mm@R=2mm{
{7}\ar@{-}[d]\ar@{-}[dr]\ar@{-}[drr] \\
{3}\ar@{-}[d]\ar@{-}[dr]\ar@{-}[drr]\ar@{-}[drrr] & {8} & {9} \\
{2}\ar@{-}[d] &  {4} & {5} & {6} \\
{1} \\
      }}
\end{equation}

Again, this algorithm can be extended to a bijection by memorizing in a second
tree of the same shape as the first one, where letter $i$ appears in $\sigma$.
By the definition of the chain, this yields an \emph{increasing} tree, that we
shall denote by $Q_2(\sigma)$.

In our example, $Q_2(\sigma)$ is
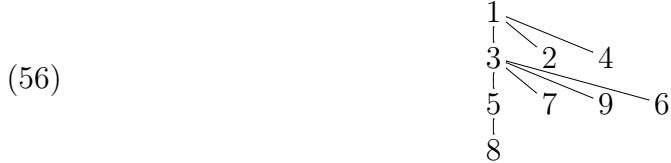
\entrymodifiers={+<4pt>}
\begin{equation}
\vcenter{\xymatrix@C=4mm@R=2mm{
{1}\ar@{-}[d]\ar@{-}[dr]\ar@{-}[drr] \\
{3}\ar@{-}[d]\ar@{-}[dr]\ar@{-}[drr]\ar@{-}[drrr] & {2} & {4} \\
{5}\ar@{-}[d] &  {7} & {9} & {6} \\
{8} \\
      }}
\end{equation}

\begin{theorem}
Consider a permutation $\sigma$.
All linear extensions of $P_2(\sigma)$, which are in bijection with the
increasing trees of the same shape as $P_2(\sigma)$, are exactly the
permutations equivalent to $\sigma$ under $\equiv_2$.
\end{theorem}

\Proof
It is obvious that two permutations having the same $\lrm$ give the same
result by the $P_2$ algorithm. Conversely, consider $P_2(\sigma)$. Given the
description of the poset, it is direct that all its linear extensions have
same $\lrm$.
\qed

Note that given a naked poset, there are two ways to label it: consider the
longest chain in it and label the last two elements either $1$ and $2$ or $2$
and $1$. This being fixed, all the other vertices have no choice  for their labeling.

\begin{corollary}[Prop. 3.11 of~\cite{Kus}]\label{cor:question}
The multisets of sizes of classes in $S_n$ under both equivalences are the
same.
\end{corollary}

\Proof
Both equivalences $\equiv$ and $\equiv_2$ give rise to the same naked posets
and all posets have, in each case, exactly two different labelings.
\qed

Finally, it is pretty clear that one could transform these equalities into a
bijection between classes, \emph{e.g.}, sending a poset of $\equiv$ with $1$
below $n$ to the poset of $\equiv_2$ of the same shape with $1$ below $2$, and
then extending this bijection to permutations by reading the same linear
extensions out of these (in other words, the one-to-one correspondence relates
permutations having same $Q$-symbols).

%%%%%%%%%%%%%%%%%%%%%%%%%%%%%%%%%%%%%%%%%%%%%%%%%%%%%%%%%%%%%%%%%%%%%%%%%%%%%%%%
\subsection{The $\{\{321,231\},\{312,132\}\}$-equivalence}

These relations are the mirror-image of those defining $\equiv_2$, which
transforms the equivalence relation on the inverse permutations into that of
having the same left-to right maxima, that is, the same saillance
composition, whence Proposition \ref{prop:congr}.

%%%%%%%%%%%%%%%%%%%%%%%%%%%%%%%%%%%%%%%%%%%%%%%%%%%%%%%%%%%%%%%%%%%%%%%%%%%%%%%
%%%%%%%%%%%%%%%%%%%%%%%%%%%%%%%%%%%%%%%%%%%%%%%%%%%%%%%%%%%%%%%%%%%%%%%%%%%%%%%

\newpage

%%%%%%%%%%%%%%%%%%%%%%%%%%%%%%%%%%%%%%%%%%%%%%%%%%%%%%%%%%%%%%%%%%%%%%%%%%%%%%%%
%%%%%%%%%%%%%%%%%%%%%%%%%%%%%%%%%%%%%%%%%%%%%%%%%%%%%%%%%%%%%%%%%%%%%%%%%%%%%%%%
%%%%%%%%%%%%%%%%%%%%%%%%%%%%%%%%%%%%%%%%%%%%%%%%%%%%%%%%%%%%%%%%%%%%%%%%%%%%%%%%
\section{Tables}

\begin{figure}[ht]
\footnotesize
%\begin{mdframed}
\begin{align*}
U_{2} &=  F_{2}\\
U_{1  1} &=  F_{1  1}\\
\\
U_{3} &=  F_{3}\\
U_{2  1} &=  F_{1  2} + F_{2  1}\\
U_{1  2} &=  F_{1  2}\\
U_{1  1  1} &=  F_{1  1  1} + F_{2  1}\\
\\
U_{4} &=  F_{4}\\
U_{3  1} &=  F_{1  3} + F_{2  2} + F_{3  1}\\
U_{2  2} &=  F_{1  3} + F_{2  2}\\
U_{2  1  1} &=  F_{1  1  2} + 2 F_{1  2  1} + F_{2  1  1} + F_{2  2} + F_{3  1}\\
U_{1  3} &=  F_{1  3}\\
U_{1  2  1} &=  F_{1  1  2} + F_{1  2  1} + F_{2  2}\\
U_{1  1  2} &=  F_{1  1  2} + F_{2  2}\\
U_{1  1  1  1} &=  F_{1  1  1  1} + 2 F_{1  2  1} + 2 F_{2  1  1} + F_{3  1}\\
\\
U_{5} &=  F_{5}\\
U_{4  1} &=  F_{1  4} + F_{2  3} + F_{3  2} + F_{4  1}\\
U_{3  2} &=  F_{1  4} + F_{2  3} + F_{3  2}\\
U_{3  1  1} &=  F_{1  1  3} + 2 F_{1  2  2} + 2 F_{1  3  1} + F_{2  1  2} + 2 F_{2  2  1} + F_{2  3} + F_{3  1  1} + F_{3  2} + F_{4  1}\\
U_{2  3} &=  F_{1  4} + F_{2  3}\\
U_{2  2  1} &=  F_{1  1  3} + 2 F_{1  2  2} + F_{1  3  1} + F_{2  1  2} + F_{2  2  1} + F_{2  3} + F_{3  2}\\
U_{2  1  2} &=  F_{1  1  3} + 2 F_{1  2  2} + F_{2  1  2} + F_{2  3} + F_{3  2}\\
U_{2  1  1  1} &=  F_{1  1  1  2} + 3 F_{1  1  2  1} + 3 F_{1  2  1  1} + 2 F_{1  2  2} + 3 F_{1  3  1} + F_{2  1  1  1} + 2 F_{2  1  2} + 5 F_{2  2  1} + 2 F_{3  1  1} + F_{3  2} + F_{4  1}\\
U_{1  4} &=  F_{1  4}\\
U_{1  3  1} &=  F_{1  1  3} + F_{1  2  2} + F_{1  3  1} + F_{2  3}\\
U_{1  2  2} &=  F_{1  1  3} + F_{1  2  2} + F_{2  3}\\
U_{1  2  1  1} &=  F_{1  1  1  2} + 2 F_{1  1  2  1} + F_{1  2  1  1} + 2 F_{1  2  2} + F_{1  3  1} + 2 F_{2  1  2} + 2 F_{2  2  1} + F_{3  2}\\
U_{1  1  3} &=  F_{1  1  3} + F_{2  3}\\
U_{1  1  2  1} &=  F_{1  1  1  2} + F_{1  1  2  1} + 2 F_{1  2  2} + 2 F_{2  1  2} + F_{2  2  1} + F_{3  2}\\
U_{1  1  1  2} &=  F_{1  1  1  2} + 2 F_{1  2  2} + 2 F_{2  1  2} + F_{3  2}\\
U_{1  1  1  1  1} &=  F_{1  1  1  1  1} + 3 F_{1  1  2  1} + 5 F_{1  2  1  1} + 3 F_{1  3  1} + 3 F_{2  1  1  1} + 5 F_{2  2  1} + 3 F_{3  1  1} + F_{4  1}\\
\end{align*}
%\end{mdframed}
\caption{The basis $U_I$ for $n\le 5$.}\label{fig:baseU}
\end{figure}
%\newpage

\def\myhs{\hskip4pt}
\def\myhsd{\hskip1.5mm}
\def\myhst{\hskip2.3mm}

\begin{figure}[ht]
\footnotesize
%	\begin{mdframed}
$$
\begin{array}{p{4mm}p{4mm}p{4mm}p{4mm}}%rrrr}
3 & 12 & 21 & 111 \\
%,\ 12,\ 21,\ 1^3
\end{array}
$$ 
$$
\left(\begin{array}{p{4mm}p{4mm}p{4mm}p{4mm}}%rrrr}
1 & 0 & 0 & 0 \\
0 & 1 & 1 & 0 \\
0 & 0 & 1 & 1 \\
0 & 0 & 0 & 1
\end{array}\right)
$$
\bigskip
$$
\begin{array}{p{5mm}p{5mm}p{5mm}p{5mm}p{5mm}p{5mm}p{5mm}p{5mm}}
4 & 13 & 22 & 31 & 112 &121 &211 &1111
\end{array}
$$
$$
\left(\begin{array}{p{5mm}p{5mm}p{5mm}p{5mm}p{5mm}p{5mm}p{5mm}p{5mm}}
1 & 0 & 0 & 0 & 0 & 0 & 0 & 0 \\
0 & 1 & 1 & 1 & 0 & 0 & 0 & 0 \\
0 & 0 & 1 & 1 & 1 & 1 & 1 & 0 \\
0 & 0 & 0 & 1 & 0 & 0 & 1 & 1 \\
0 & 0 & 0 & 0 & 1 & 1 & 1 & 0 \\
0 & 0 & 0 & 0 & 0 & 1 & 2 & 2 \\
0 & 0 & 0 & 0 & 0 & 0 & 1 & 2 \\
0 & 0 & 0 & 0 & 0 & 0 & 0 & 1
\end{array}\right)
$$
\bigskip
$$
\begin{array}{p{5mm}p{5mm}p{5mm}p{5mm}p{5mm}p{5mm}p{5mm}p{5mm}
p{5mm}p{5mm}p{5mm}p{5mm}p{5mm}p{5mm}p{5mm}p{5mm}}
5 & 14 & 23 & 32 & 41 & 113 & 122 & 131 & 212 & 221 & 311 & 1112
& 1121 & 1211 & 2111 & 11111
\end{array}
$$ 
$$
\left(\begin{array}{p{5mm}p{5mm}p{5mm}p{5mm}p{5mm}p{5mm}p{5mm}p{5mm}
p{5mm}p{5mm}p{5mm}p{5mm}p{5mm}p{5mm}p{5mm}p{5mm}}
%\left(\begin{array}{rrrrrrrrrrrrrrrr}
1 & 0 & 0 & 0 & 0 & 0 & 0 & 0 & 0 & 0 & 0 & 0 & 0 & 0 & 0 & 0\\\
0 & 1 & 1 & 1 & 1 & 0 & 0 & 0 & 0 & 0 & 0 & 0 & 0 & 0 & 0 & 0 \\
0 & 0 & 1 & 1 & 1 & 1 & 1 & 1 & 1 & 1 & 1 & 0 & 0 & 0 & 0 & 0 \\
0 & 0 & 0 & 1 & 1 & 0 & 0 & 0 & 1 & 1 & 1 & 1 & 1 & 1 & 1 & 0 \\
0 & 0 & 0 & 0 & 1 & 0 & 0 & 0 & 0 & 0 & 1 & 0 & 0 & 0 & 1 & 1 \\
0 & 0 & 0 & 0 & 0 & 1 & 1 & 1 & 1 & 1 & 1 & 0 & 0 & 0 & 0 & 0 \\
0 & 0 & 0 & 0 & 0 & 0 & 1 & 1 & 2 & 2 & 2 & 2 & 2 & 2 & 2 & 0 \\
0 & 0 & 0 & 0 & 0 & 0 & 0 & 1 & 0 & 1 & 2 & 0 & 0 & 1 & 3 & 3 \\
0 & 0 & 0 & 0 & 0 & 0 & 0 & 0 & 1 & 1 & 1 & 2 & 2 & 2 & 2 & 0 \\
0 & 0 & 0 & 0 & 0 & 0 & 0 & 0 & 0 & 1 & 2 & 0 & 1 & 2 & 5 & 5 \\
0 & 0 & 0 & 0 & 0 & 0 & 0 & 0 & 0 & 0 & 1 & 0 & 0 & 0 & 2 & 3 \\
0 & 0 & 0 & 0 & 0 & 0 & 0 & 0 & 0 & 0 & 0 & 1 & 1 & 1 & 1 & 0 \\
0 & 0 & 0 & 0 & 0 & 0 & 0 & 0 & 0 & 0 & 0 & 0 & 1 & 2 & 3 & 3 \\
0 & 0 & 0 & 0 & 0 & 0 & 0 & 0 & 0 & 0 & 0 & 0 & 0 & 1 & 3 & 5 \\
0 & 0 & 0 & 0 & 0 & 0 & 0 & 0 & 0 & 0 & 0 & 0 & 0 & 0 & 1 & 3 \\
0 & 0 & 0 & 0 & 0 & 0 & 0 & 0 & 0 & 0 & 0 & 0 & 0 & 0 & 0 & 1
\end{array}\right)
$$
%	\end{mdframed}
	\caption{Transition matrices $M_n$ for $n\le 5$.}\label{fig:M}
\end{figure} 
%\newpage

\begin{figure}[ht]
\footnotesize
%	\begin{mdframed}
\begin{align*}
V_{2} &=  R_{2}\\
V_{1  1} &=  R_{1  1}\\
\\
\\
V_{3} &=  R_{3}\\
V_{2  1} &=  -R_{1  1  1} + R_{2  1}\\
V_{1  2} &=  R_{1  1  1} + R_{1  2} - R_{2  1}\\
V_{1  1  1} &=  R_{1  1  1}\\
\\
\\
V_{4} &=  R_{4}\\
V_{3  1} &=  R_{1  1  1  1} - R_{2  1  1} + R_{3  1}\\
V_{2  2} &=  -R_{1  1  1  1} - R_{1  1  2} + R_{2  1  1} + R_{2  2} - R_{3  1}\\
V_{2  1  1} &=  -2 R_{1  1  1  1} + R_{2  1  1}\\
V_{1  3} &=  R_{1  1  2} + R_{1  3} - R_{2  2}\\
V_{1  2  1} &=  2 R_{1  1  1  1} + R_{1  2  1} - 2 R_{2  1  1}\\
V_{1  1  2} &=  R_{1  1  2} - R_{1  2  1} + R_{2  1  1}\\
V_{1  1  1  1} &=  R_{1  1  1  1}\\
\\
\\
V_{5} &=  R_{5}\\
V_{4  1} &=  -R_{1  1  1  1  1} + R_{2  1  1  1} - R_{3  1  1} + R_{4  1}\\
V_{3  2} &=  R_{1  1  1  1  1} + R_{1  1  1  2} - R_{2  1  1  1} - R_{2  1  2} + R_{3  1  1} + R_{3  2} - R_{4  1}\\
V_{3  1  1} &=  3 R_{1  1  1  1  1} - 2 R_{2  1  1  1} + R_{3  1  1}\\
V_{2  3} &=  -R_{1  1  1  2} - R_{1  1  3} + R_{2  1  2} + R_{2  3} - R_{3  2}\\
V_{2  2  1} &=  -2 R_{1  1  1  1  1} - R_{1  1  2  1} + 2 R_{2  1  1  1} + R_{2  2  1} - 2 R_{3  1  1}\\
V_{2  1  2} &=  -R_{1  1  1  1  1} - 2 R_{1  1  1  2} + R_{1  1  2  1} + R_{2  1  2} - R_{2  2  1} + R_{3  1  1}\\
V_{2  1  1  1} &=  -3 R_{1  1  1  1  1} + R_{2  1  1  1}\\
V_{1  4} &=  R_{1  1  3} + R_{1  4} - R_{2  3}\\
V_{1  3  1} &=  -2 R_{1  1  1  1  1} + R_{1  1  2  1} - R_{1  2  1  1} + R_{1  3  1} + 2 R_{2  1  1  1} - R_{2  2  1}\\
V_{1  2  2} &=  2 R_{1  1  1  1  1} + 2 R_{1  1  1  2} - R_{1  1  2  1} + R_{1  2  1  1} + R_{1  2  2} - R_{1  3  1} - 2 R_{2  1  1  1} - 2 R_{2  1  2} + R_{2  2  1}\\
V_{1  2  1  1} &=  4 R_{1  1  1  1  1} + R_{1  2  1  1} - 3 R_{2  1  1  1}\\
V_{1  1  3} &=  R_{1  1  3} - R_{1  2  2} + R_{2  1  2}\\
V_{1  1  2  1} &=  -2 R_{1  1  1  1  1} + R_{1  1  2  1} - 2 R_{1  2  1  1} + 3 R_{2  1  1  1}\\
V_{1  1  1  2} &=  R_{1  1  1  1  1} + R_{1  1  1  2} - R_{1  1  2  1} + R_{1  2  1  1} - R_{2  1  1  1}\\
V_{1  1  1  1  1} &=  R_{1  1  1  1  1}\\
\end{align*}
%	\end{mdframed}
	\caption{The basis $V_I$ for $n\le 5$.}\label{fig:VR} 
\end{figure}
\pagebreak[3]

\begin{figure}[ht]
%	\begin{mdframed}
\footnotesize
\begin{align*}
V_{2} &=  \Lambda_{1  1} - \Lambda_{2}\\
V_{1  1} &=  \Lambda_{2}\\
\\
\\
V_{3} &=  \Lambda_{1  1  1} - \Lambda_{1  2} - \Lambda_{2  1} + \Lambda_{3}\\
V_{2  1} &=  \Lambda_{1  2} - 2 \Lambda_{3}\\
V_{1  2} &=  -\Lambda_{1  2} + \Lambda_{2  1} + \Lambda_{3}\\
V_{1  1  1} &=  \Lambda_{3}\\
\\
\\
V_{4} &=  \Lambda_{1  1  1  1} - \Lambda_{1  1  2} - \Lambda_{1  2  1} + \Lambda_{1  3} - \Lambda_{2  1  1} + \Lambda_{2  2} + \Lambda_{3  1} - \Lambda_{4}\\
V_{3  1} &=  \Lambda_{1  1  2} - 2 \Lambda_{1  3} - \Lambda_{2  2} + 3 \Lambda_{4}\\
V_{2  2} &=  -\Lambda_{1  1  2} + \Lambda_{1  2  1} + \Lambda_{1  3} + \Lambda_{2  2} - 2 \Lambda_{3  1} - \Lambda_{4}\\
V_{2  1  1} &=  \Lambda_{1  3} - 3 \Lambda_{4}\\
V_{1  3} &=  -\Lambda_{1  2  1} + \Lambda_{1  3} + \Lambda_{2  1  1} - \Lambda_{2  2} + \Lambda_{3  1} - \Lambda_{4}\\
V_{1  2  1} &=  -2 \Lambda_{1  3} + \Lambda_{2  2} + 3 \Lambda_{4}\\
V_{1  1  2} &=  \Lambda_{1  3} - \Lambda_{2  2} + \Lambda_{3  1} - \Lambda_{4}\\
V_{1  1  1  1} &=  \Lambda_{4}\\
\\
\\
V_{5} &=  \Lambda_{1  1  1  1  1} - \Lambda_{1  1  1  2} - \Lambda_{1  1  2  1} + \Lambda_{1  1  3} - \Lambda_{1  2  1  1} + \Lambda_{1  2  2}\\
	&+ \Lambda_{1  3  1} - \Lambda_{1  4} - \Lambda_{2  1  1  1} + \Lambda_{2  1  2} + \Lambda_{2  2  1} - \Lambda_{2  3} + \Lambda_{3  1  1} - \Lambda_{3  2} - \Lambda_{4  1} + \Lambda_{5}\\
V_{4  1} &=  \Lambda_{1  1  1  2} - 2 \Lambda_{1  1  3} - \Lambda_{1  2  2} + 3 \Lambda_{1  4} - \Lambda_{2  1  2} + 2 \Lambda_{2  3} + \Lambda_{3  2} - 4 \Lambda_{5}\\
V_{3  2} &=  -\Lambda_{1  1  1  2} + \Lambda_{1  1  2  1} + \Lambda_{1  1  3} + \Lambda_{1  2  2} - 2 \Lambda_{1  3  1} - \Lambda_{1  4} + \Lambda_{2  1  2} - \Lambda_{2  2  1} - \Lambda_{2  3} - \Lambda_{3  2} + 3 \Lambda_{4  1} + \Lambda_{5}\\
V_{3  1  1} &=  \Lambda_{1  1  3} - 3 \Lambda_{1  4} - \Lambda_{2  3} + 6 \Lambda_{5}\\
V_{2  3} &=  -\Lambda_{1  1  2  1} + \Lambda_{1  1  3} + \Lambda_{1  2  1  1} - \Lambda_{1  2  2} + \Lambda_{1  3  1} - \Lambda_{1  4} + \Lambda_{2  2  1} - \Lambda_{2  3} - 2 \Lambda_{3  1  1} + 2 \Lambda_{3  2} - \Lambda_{4  1} + \Lambda_{5}\\
V_{2  2  1} &=  -2 \Lambda_{1  1  3} + \Lambda_{1  2  2} + 3 \Lambda_{1  4} + 2 \Lambda_{2  3} - 2 \Lambda_{3  2} - 4 \Lambda_{5}\\
V_{2  1  2} &=  \Lambda_{1  1  3} - \Lambda_{1  2  2} + \Lambda_{1  3  1} - \Lambda_{1  4} - \Lambda_{2  3} + 2 \Lambda_{3  2} - 3 \Lambda_{4  1} + \Lambda_{5}\\
V_{2  1  1  1} &=  \Lambda_{1  4} - 4 \Lambda_{5}\\
V_{1  4} &=  -\Lambda_{1  2  1  1} + \Lambda_{1  2  2} + \Lambda_{1  3  1} - \Lambda_{1  4} + \Lambda_{2  1  1  1} - \Lambda_{2  1  2} - \Lambda_{2  2  1} + \Lambda_{2  3} + \Lambda_{3  1  1} - \Lambda_{3  2} - \Lambda_{4  1} + \Lambda_{5}\\
V_{1  3  1} &=  -\Lambda_{1  2  2} + 3 \Lambda_{1  4} + \Lambda_{2  1  2} - 2 \Lambda_{2  3} + \Lambda_{3  2} - 4 \Lambda_{5}\\
V_{1  2  2} &=  \Lambda_{1  2  2} - 2 \Lambda_{1  3  1} - \Lambda_{1  4} - \Lambda_{2  1  2} + \Lambda_{2  2  1} + \Lambda_{2  3} - \Lambda_{3  2} + 3 \Lambda_{4  1} + \Lambda_{5}\\
V_{1  2  1  1} &=  -3 \Lambda_{1  4} + \Lambda_{2  3} + 6 \Lambda_{5}\\
V_{1  1  3} &=  \Lambda_{1  3  1} - \Lambda_{1  4} - \Lambda_{2  2  1} + \Lambda_{2  3} + \Lambda_{3  1  1} - \Lambda_{3  2} - \Lambda_{4  1} + \Lambda_{5}\\
V_{1  1  2  1} &=  3 \Lambda_{1  4} - 2 \Lambda_{2  3} + \Lambda_{3  2} - 4 \Lambda_{5}\\
V_{1  1  1  2} &=  -\Lambda_{1  4} + \Lambda_{2  3} - \Lambda_{3  2} + \Lambda_{4  1} + \Lambda_{5}\\
V_{1  1  1  1  1} &=  \Lambda_{5}\\
\end{align*}
%	\end{mdframed}
\caption{The basis $V_I$ on the $\Lambda^J$ for $n\le 5$.}\label{fig:VL} 
\end{figure}
%\newpage

\begin{figure}[ht]
%	\begin{mdframed}
\footnotesize
$$
2\ \ 11
$$
$$
\left(\begin{array}{rr}
1 & 0 \\
0 & 1
\end{array}\right)
$$

\bigskip

$$
\begin{array}{p{4mm}p{4mm}p{4mm}p{4mm}}%rrrr}
3 & 12 & 21 & 111 \\
%,\ 12,\ 21,\ 1^3
\end{array}
$$
$$
\left(\begin{array}{rrrr}
1 & 0 & 0 & 0 \\
0 & 1 & -1 & 1 \\
0 & 0 & 1 & -1 \\
0 & 0 & 0 & 1
\end{array}\right)
$$

\bigskip

$$
\begin{array}{p{5mm}p{5mm}p{5mm}p{5mm}p{5mm}p{5mm}p{5mm}p{5mm}}
4 & 13 & 22 & 31 & 112 &121 &211 &1111
\end{array}
$$
$$
\left(\begin{array}{p{5mm}p{5mm}p{5mm}p{5mm}p{5mm}p{5mm}p{5mm}p{5mm}}
1 & 0 & 0 & 0 & 0 & 0 & 0 & 0 \\
0 & 1 & -1 & 0 & 1 & 0 & 0 & 0 \\
0 & 0 & 1 & -1 & -1 & 0 & 1 & -1 \\
0 & 0 & 0 & 1 & 0 & 0 & -1 & 1 \\
0 & 0 & 0 & 0 & 1 & -1 & 1 & 0 \\
0 & 0 & 0 & 0 & 0 & 1 & -2 & 2 \\
0 & 0 & 0 & 0 & 0 & 0 & 1 & -2 \\
0 & 0 & 0 & 0 & 0 & 0 & 0 & 1
\end{array}\right)
$$
\bigskip
$$
\begin{array}{p{5mm}p{5mm}p{5mm}p{5mm}p{5mm}p{5mm}p{5mm}p{5mm}
p{5mm}p{5mm}p{5mm}p{5mm}p{5mm}p{5mm}p{5mm}p{5mm}}
5 & 14 & 23 & 32 & 41 & 113 & 122 & 131 & 212 & 221 & 311 & 1112
& 1121 & 1211 & 2111 & 11111
\end{array}
$$
$$
\left(\begin{array}{p{5mm}p{5mm}p{5mm}p{5mm}p{5mm}p{5mm}p{5mm}p{5mm}
p{5mm}p{5mm}p{5mm}p{5mm}p{5mm}p{5mm}p{5mm}p{5mm}}
%\left(\begin{array}{rrrrrrrrrrrrrrrr}
1 & 0 & 0 & 0 & 0 & 0 & 0 & 0 & 0 & 0 & 0 & 0 & 0 & 0 & 0 & 0 \\
0 & 1 & -1 & 0 & 0 & 1 & 0 & 0 & 0 & 0 & 0 & 0 & 0 & 0 & 0 & 0 \\
0 & 0 & 1 & -1 & 0 & -1 & 0 & 0 & 1 & 0 & 0 & -1 & 0 & 0 & 0 & 0 \\
0 & 0 & 0 & 1 & -1 & 0 & 0 & 0 & -1 & 0 & 1 & 1 & 0 & 0 & -1 & 1 \\
0 & 0 & 0 & 0 & 1 & 0 & 0 & 0 & 0 & 0 & -1 & 0 & 0 & 0 & 1 & -1 \\
0 & 0 & 0 & 0 & 0 & 1 & -1 & 0 & 1 & 0 & 0 & 0 & 0 & 0 & 0 & 0 \\
0 & 0 & 0 & 0 & 0 & 0 & 1 & -1 & -2 & 1 & 0 & 2 & -1 & 1 & -2 & 2 \\
0 & 0 & 0 & 0 & 0 & 0 & 0 & 1 & 0 & -1 & 0 & 0 & 1 & -1 & 2 & -2 \\
0 & 0 & 0 & 0 & 0 & 0 & 0 & 0 & 1 & -1 & 1 & -2 & 1 & 0 & 0 & -1 \\
0 & 0 & 0 & 0 & 0 & 0 & 0 & 0 & 0 & 1 & -2 & 0 & -1 & 0 & 2 & -2 \\
0 & 0 & 0 & 0 & 0 & 0 & 0 & 0 & 0 & 0 & 1 & 0 & 0 & 0 & -2 & 3 \\
0 & 0 & 0 & 0 & 0 & 0 & 0 & 0 & 0 & 0 & 0 & 1 & -1 & 1 & -1 & 1 \\
0 & 0 & 0 & 0 & 0 & 0 & 0 & 0 & 0 & 0 & 0 & 0 & 1 & -2 & 3 & -2 \\
0 & 0 & 0 & 0 & 0 & 0 & 0 & 0 & 0 & 0 & 0 & 0 & 0 & 1 & -3 & 4 \\
0 & 0 & 0 & 0 & 0 & 0 & 0 & 0 & 0 & 0 & 0 & 0 & 0 & 0 & 1 & -3 \\
0 & 0 & 0 & 0 & 0 & 0 & 0 & 0 & 0 & 0 & 0 & 0 & 0 & 0 & 0 & 1
\end{array}\right)
$$
%\end{mdframed}
\caption{The matrices $M_n^{-1}$ for $n\le 5$.}
\end{figure}

\end{document}